\documentclass[letterpaper,12pt]{article}

\usepackage{color}
\usepackage{amssymb}
\usepackage{amsmath}
\usepackage{amscd}
\usepackage{mathtools}
\usepackage[all]{xy}
\usepackage{tikz-cd}

\usepackage{hyperref}
\hypersetup{colorlinks=true,linkcolor=blue,anchorcolor=blue,citecolor=blue}
\usepackage{cleveref}

\usepackage{color}

\newtheorem{theo}{Theorem}[section]
\newtheorem{prop}[theo]{Proposition}
\newtheorem{lem}[theo]{Lemma}
\newtheorem{cor}[theo]{Corollary}
\newtheorem{defi}[theo]{Definition}
\newtheorem{rem}[theo]{Remark}
\newtheorem{conj}[theo]{Conjecture}

\def\id{{\rm id}}
\def \fppf {{\rm fppf}}

\def \cris {{\rm cris}}
\def \s {{\rm s}}

\def \Br {{\rm{Br}}}

\def \Ga {{\Gamma}}

\def \Pic {{\rm {Pic}}}

\def \Gal {{\rm{Gal}}}
\def \Ker {{\rm{Ker}}}
\def \Im {{\rm {Im}}}

\def \P{{\mathbb P}}
\def \Spec {{\rm{Spec}}}
\def \dim {{\rm{dim}}}
\def \Hom {{\rm {Hom}}}
\def \End {{\rm {End}}}

\def \Z {{\mathbb Z}}
\def \Q {{\mathbb Q}}
\def \F {{\mathbb F}}
\def \Ext {{\rm Ext}}

\def\G{{\mathbb G}}

\def\GG{{\cal G}}

\def\T{{\cal T}}

\def\lra{\longrightarrow}

\def\H{{\rm H}}

\def\Zar{{\rm Zar}}

\def\Kum{{\rm Kum}}

\def\NS{{\rm NS\,}}

\def\O{{\cal O}}

\def\sE{{\cal E}}

\def\res{{\rm res}}

\def\Ga{\Gamma}

\def\e{\varepsilon}
\def\et{{\rm{\acute et}}}

\def\prim{{\rm prim}}

\newcommand{\bthe}{\begin{theo}}
\newcommand{\ble}{\begin{lem}}
\newcommand{\bpr}{\begin{prop}}
\newcommand{\bco}{\begin{cor}}
\newcommand{\bde}{\begin{defi}}
\newcommand{\brem}{\begin{rem}}
\newcommand{\bconj}{\begin{conj}}
\newcommand{\ethe}{\end{theo}}
\newcommand{\ele}{\end{lem}}
\newcommand{\epr}{\end{prop}}
\newcommand{\eco}{\end{cor}}
\newcommand{\ede}{\end{defi}}
\newcommand{\erem}{\end{rem}}
\newcommand{\econj}{\end{conj}}

\newcommand\blfootnote[1]{%
  \begingroup
  \renewcommand\thefootnote{}\footnote{#1}%
  \addtocounter{footnote}{-1}%
  \endgroup
}

\title{Boundedness of the $p$-primary torsion of the Brauer group of products of varieties}
\author{Alexei N.~Skorobogatov}

\date{\today}
\begin{document}
\maketitle

\blfootnote{The 2020 Subject Classification codes: 14F22, 14K15, 14F30}

\begin{abstract}
\noindent Let $k$ be a field finitely generated over its prime subfield. 
We prove that the quotient of the Brauer group of a product of varieties over $k$ by 
the sum of the images of the Brauer
groups of factors has finite exponent. The bulk of the proof concerns $p$-primary torsion
in characteristic $p$. Our approach gives a more direct proof of 
the boundedness of the $p$-primary torsion of the Brauer group of an abelian variety, as recently
proved by D'Addezio. We show that
the transcendental Brauer group of a Kummer surface over $k$ has finite exponent, but
can be infinite when $k$ is an infinite field of positive characteristic. 
This answers a question of Zarhin and the author. 
\end{abstract}

\section*{Introduction}

Let $k$ be a field of characteristic exponent $p$.
Thus $p=1$ if ${\rm char}(k)=0$, otherwise $p={\rm char}(k)$.
Let $\bar k$ be an algebraic closure of $k$, 
let $k^\s$ be the separable closure of $k$ in $\bar k$, and let 
$\Ga=\Gal(k^\s/k)$. For an abelian group $A$ and a prime number
$\ell$ we denote by $A\{\ell\}$ the $\ell$-primary torsion subgroup of $A$.
We write $A(p')$ for the direct sum of $A\{\ell\}$ over all primes $\ell\neq p$.

Assume that $k$ is finitely generated over its prime subfield.
Relation between the Tate conjecture for divisors 
for a smooth and projective variety $X$ over $k$ and finiteness properties of the Brauer group of $X$
is well known, at least for torsion coprime to $p$. Indeed, the validity
of the Tate conjecture for $X$ at a prime $\ell\neq p$ 
is equivalent to the finiteness of $\Br(X_{k^\s})^\Ga\{\ell\}$, and is also equivalent 
to the finiteness of the image of the natural map $\Br(X)\{\ell\}\to\Br(X_{k^\s})\{\ell\}$,
see \cite[Thm.~16.1.1]{CTS21}. In particular, this holds for abelian varieties and K3 surfaces.
Moreover, in these two cases $\Br(X_{k^\s})^\Ga(p')$ 
is finite \cite{SZ08, SZ15, Ito}, see also \cite[Ch.~16]{CTS21}.
In  \cite[Questions 1, 2]{SZ08} the authors asked whether $\Br(X_{k^\s})^\Ga\{p\}$, or at least 
the image of $\Br(X)\{p\}$ in $\Br(X_{k^\s})\{p\}$, is finite when $X$ is an abelian variety or
a K3 surface and $p>1$. In a recent paper, D'Addezio
observed that for the self-product of a supersingular elliptic curve
this image is infinite when $k$ is infinite \cite[Cor.~5.4]{D'A}. On the positive side, he proved that 
$\Br(X_{k^\s})^\Ga\{p\}$ has {\em finite exponent} when $X$ is an abelian variety, see \cite[Thm.~1.1]{D'A}. 
(As pointed out in \cite[Cor.~6.7]{D'A}, this may fail if $k^\s$ is replaced by $\bar k$.)
For $p\neq 2$, we note that D'Addezio's examples descend to the associated Kummer surfaces.
Thus the questions raised in \cite{SZ08} have negative answers for K3 surfaces 
over infinite finitely generated fields of characteristic $p\geq 3$.

The main result of this note is the following

\medskip

\noindent{\bf Theorem A}
{\em Let $X$ and $Y$ be smooth, projective, geometrically integral
varieties over a finitely generated field $k$. Then
the cokernel of the natural map $$\Br(X)\oplus\Br(Y)\to\Br(X\times_kY)$$
has finite exponent.}

\medskip

For the prime-to-$p$ torsion this easily follows from
\cite[Thm.~B]{SZ14} which says\footnote{In {\em loc.~cit.}~one assumes that ${\rm char}(k)=0$, but the proof goes through for the prime-to-$p$
torsion when ${\rm char}(k)=p>0$.} that the cokernel of 
$\Br(X)(p')\oplus\Br(Y)(p')\to\Br(X\times_kY)(p')$ is finite when $X\times_kY$ has a $k$-point
or $\H^3(k,(k^\s)^\times)=0$. In this paper we deal with the $p$-primary torsion.
Our proof is inspired by \cite{D'A} and crucially uses
the crystalline Tate conjecture proved by de Jong \cite[Thm.~2.6]{dJ98}. 
As a consequence we obtain a more transparent proof of \cite[Thm.~1.1]{D'A}.
Combined with the previous results of Zarhin and the author, it gives that $\Br(X_{k^\s})^\Ga$
 is a direct sum of a finite group and a $p$-group of finite exponent 
when $X$ is an abelian variety over a finitely generated field $k$, see Theorem \ref{da}.  
Using similar ideas, we also give a simplified proof of the flat version of the Tate conjecture
for divisors on abelian varieties \cite[Thm.~5.1]{D'A}, see Theorem \ref{dada}.

The prime-to-$p$ torsion part of the next result was obtained in \cite[Thm.~3.1]{SZ14}.

\medskip

\noindent{\bf Theorem B}
{\em Let $X$ and $Y$ be smooth, projective, geometrically integral
varieties over a finitely generated field $k$ of characteristic exponent $p$.
Then the cokernel of the natural map
$\Br(X_{k^\s})^\Ga\oplus\Br(Y_{k^\s})^\Ga\to\Br(X_{k^\s}\times_{k^\s}Y_{k^\s})^\Ga$ 
is a direct sum of a finite group and a $p$-group of finite exponent.}

\medskip

Theorem B can be used to prove that for some surfaces $X$ 
dominated by a product of curves, $\Br(X_{k^\s})^\Ga$
is a direct sum of a finite group and a $p$-group of finite exponent, see Corollary \ref{chennai}.

Our approach is based on the systematic use of {\em pointed varieties}, i.e.~varieties over
$k$ with a distinguished $k$-point. In Section \ref{1} we obtain
a version of the K\"unneth formula for the second flat cohomology group
of the product of pointed varieties, see Theorem \ref{ku} and Corollary \ref{si}.
Similarly to the $\ell$-adic case, the embedding of the `primitive' part of cohomology
can be interpreted in terms of pairing with classes of certain natural torsors.
In Section \ref{3} we first prove Theorem A for pointed varieties (Theorem \ref{kanchi})
from which we obtain the general case, see Theorem \ref{wed}. We then deduce
Theorem B, see Corollary \ref{wed1}.
Applications to abelian varieties can be found in Section \ref{A} and applications
to Kummer surfaces in Section \ref{K}. We show that
the transcendental Brauer group of a Kummer surface over a finitely generated field
$k$ of characteristic not equal to 2 has finite exponent, but is infinite for the Kummer surface
attached to the self-product of a supersingular elliptic curve when
$k$ is infinite of positive characteristic. 

The appendix by Alexander Petrov
contains a structure theorem for the $p$-primary torsion subgroup of
the Brauer group of a smooth and proper variety
over an algebraically closed field of positive characteristic $p$: this group is a direct sum
of finitely many copies of $\Q_p/\Z_p$ and an abelian $p$-group of finite exponent,
see Theorem \ref{appendix: brauer finite}. 
This can be deduced from \cite[Prop.~II.5.9]{Ill} and its proof,
and seems to be well known to the experts.
The proof of Theorem \ref{appendix: brauer finite} given
in the appendix is a self-contained argument that relies on some
basic properties of the de Rham--Witt complex.

The work on this paper started when the author visited Capital Normal University 
in Beijing and continued during visits to Chennai Mathematical Institute and EPF Lausanne. 
He is grateful to Yang Cao, Marco D'Addezio, Jean-Pierre Serre, Domenico Valloni, 
Yuan Yang, and Yuri Zarhin for stimulating discussions, and to Alexander Petrov who very kindly
provided the appendix to this paper.

\section{Cohomology of the product} \label{1}

Let $k$ be a field. 
Let $F$ be a contravariant functor from the category of 
schemes over $k$ to the category of abelian groups. 
We shall refer to a pair $(X,x_0)$, where 
$X$ is a $k$-scheme and 
$x_0\in X(k)$, as a {\em pointed $k$-scheme}. 
For a pointed $k$-scheme $(X,x_0)$ we
define $$F(X)_e:=\Ker[x_0^*:F(X)\to F(k)].$$
Then we have $F(X)\cong F(k)\oplus F(X)_e$.
For $k$-schemes $X$ and $Y$ we have an obvious commutative diagram
$$\xymatrix{Y\ar[d]_{\pi_Y}&\ar[l]_{p_X}X\times_k Y\ar[d]^{p_Y}\\
\Spec(k)&\ar[l]_{\ \ \ \ \ \pi_X} X}$$
When $(X,x_0)$ and $(Y,y_0)$ are pointed $k$-schemes,
the $k$-points $x_0$ and $y_0$ give rise to sections to the four morphisms in this diagram.
Thus $F(k)$, $F(X)$, $F(Y)$ are direct summands of $F(X\times_kY)$ such that
$F(X)\cap F(Y)=F(k)$. Therefore, $F(X)_e$ and $F(Y)_e$ are direct summands of $F(X\times_kY)_e$
such that $F(X)_e\cap F(Y)_e=0$. It follows that
$F(X)_e\oplus F(Y)_e$ is a direct summand of $F(X\times_kY)_e$. Define
$$F(X\times_kY)_\prim:=\Ker[F(X\times_kY)_e\to F(X)_e\oplus F(Y)_e],$$
where the map $F(X\times_kY)_e\to F(X)_e$ is the specialisation at $y_0$ and
the map $F(X\times_kY)_e\to F(Y)_e$ is the specialisation at $x_0$.
This gives rise to a direct sum decomposition of abelian groups
\begin{equation}
F(X\times_kY)_e\cong F(X)_e\oplus F(Y)_e\oplus F(X\times_kY)_\prim, \label{dec}
\end{equation}
which is functorial with respect to morphisms of pointed $k$-schemes.

\medskip

For a field extension $K/k$ we define the functor $F(X_K)^k:=\Im[F(X)\to F(X_K)]$.
The group $\Br(X_{k^\s})^k$ is called the {\em transcendental Brauer group}.

\medskip

Recall that by a theorem of Grothendieck,
the Picard scheme ${\rm{\bf Pic}}_{X/k}$ exists when $X$ is proper over $k$, see
the references in \cite[Thm.~2.5.7]{CTS21}. The Picard variety of 
a smooth, projective, geometrically integral variety $X$ is the abelian variety
${\rm{\bf Pic}}^0_{X/k, {\rm red}}$, where ${\rm{\bf Pic}}^0_{X/k}$ is the connected 
component of 0. The Albanese variety $A$ is defined as the dual abelian variety of the Picard
variety of $X$ so that ${\rm{\bf Pic}}^0_{X/k, {\rm red}}\cong A^\vee$. 

\medskip 

From now on we assume that $X$ is a projective variety over a field $k$, and that $p$
is a prime number that may or may not be equal to the characteristic of $k$, unless explicitly stated 
otherwise. Throughout the paper
we consider fppf-cohomology, so we drop fppf from notation. We also
write $\H^i(X):=\H^i(X_\fppf,\mu_{p^n})$. 

Let $S_X$ be the finite commutative group $k$-scheme whose Cartier dual $S_X^\vee$ 
is the subgroup $k$-scheme  
${\bf Pic}_{X/k}[p^n]:=\Ker[{\bf Pic}_{X/k}\xrightarrow{p^n}{\bf Pic}_{X/k}]$.

\bpr \label{mur}
Let $X$ and $Y$ be pointed projective, geometrically reduced and 
geometrically connected varieties over a field $k$. Then there is a natural isomorphism
$$
\H^2(X\times_k Y,\mu_{p^n})_\prim\cong\H^1(X, S_Y^\vee)_e.
$$
\epr
{\em Proof.} For a proper, geometrically reduced and geometrically connected
$k$-variety $\pi_Y\colon Y\to\Spec(k)$ the natural map $\O_{\Spec(k)}\to \pi_{Y*}\O_Y$ is an isomorphism.
This implies that every $k$-morphism from $Y$ to an affine $k$-scheme must be constant.
In particular, the sheaf $\pi_{Y *}\mu_{p^n,Y}$ on $\Spec(k)_\fppf$ is $\mu_{p^n}$.
The Kummer sequence 
$$1\to \mu_{p^n}\to \G_{m,k}\xrightarrow{p^n}\G_{m,k}\to 1$$
is an exact sequence of sheaves on $\Spec(k)_\fppf$.
Using that the natural morphism $\G_{m,k}\to \pi_{Y*}\G_{m,Y}$ is an isomorphism,
we see that the group $k$-scheme $S_Y^\vee$ represents the sheaf
$R^1\pi_{Y*}\mu_{p^n}$ on $\Spec(k)_\fppf$.
By a theorem of Bragg and Olsson \cite[Cor.~1.4]{BO}, since $Y$ is projective,
there is an affine group $k$-scheme $G_n$ of finite type that represents the sheaf
$R^2 \pi_{Y *}\mu_{p^n}$ on $\Spec(k)_\fppf$.

Consider the spectral sequence attached to $p_Y\colon X\times_k Y\to X$:
$$E^{p,q}_2=\H^p(X,R^q p_{Y *}\mu_{p^n})\Rightarrow \H^{p+q}(X\times_k Y).$$
Since $(\id,y_0)$ is a section of $p_Y$, the canonical map 
$$\H^i(X)\cong\H^i(X,p_{Y *}\mu_{p^n})\to \H^i(X\times_kY)$$ 
is split injective for any $i\geq 0$. 
This implies that the differentials
on any page of this spectral sequence with target $\H^i(X)$ are zero for any $i\geq 0$.
It follows that we have an exact sequence
$$0\to \H^1(X, S_Y^\vee)\to\H^2(X\times_k Y)/\H^2(X)\to 
\H^0(X,G_n)\to \H^2(X, S_Y^\vee).$$
When $X=\Spec(k)$ there is a compatible exact sequence giving rise to the commutative diagram
$$\xymatrix{
0\ar[r]& \H^1(X, S_Y^\vee)\ar[r]&\H^2(X\times_k Y)_e/\H^2(X)_e\ar[r]&
\H^0(X,G_n)\ar[r]&\H^2(X, S_Y^\vee)\\
0\ar[r]& \H^1(k, S_Y^\vee)\ar[r]\ar@{^{(}->}[u]&\H^2(Y)_e\ar[r]\ar@{^{(}->}[u]&
\H^0(k,G_n)\ar[u]^\cong\ar[r]&\H^2(k, S_Y^\vee)\ar@{^{(}->}[u]}$$
All vertical maps are split injective, with splittings defined by the base point $x_0\in X(k)$. 
The map $\H^0(k,G_n)\to\H^0(X,G_n)$ is an isomorphism
since $X$ is proper, geometrically reduced and geometrically connected, and $G_n$ is affine. 
By diagram chase we obtain a natural isomorphism
$$\H^2(X\times_k Y)_e/\big(\H^2(X)_e\oplus\H^2(Y)_e\big)\cong \H^1(X, S_Y^\vee)_e.$$
This proves the proposition. $\Box$

\medskip

The following statement can be compared to \cite[Prop.~1.1]{HS13}.

\bpr \label{tt}
Let $X$ be a pointed projective, geometrically reduced and
geometrically connected variety over a field $k$.
For any finite commutative group $k$-scheme $\GG$ we have a functorial isomorphism
$$\tau\colon
\H^1(X,\GG)_e\stackrel{\sim}\lra\Hom_k(\GG^\vee,{\bf Pic}_{X/k}). $$
\epr
{\em Proof.} We adapt the method of proof of \cite[Thm.~1.5.1]{CS}. 

There is the following spectral sequence for the fppf topology:
$$\Ext^p_{k}(A,R^q\pi_{X*}B)\Rightarrow \Ext^{p+q}_{X}(\pi_X^*A,B),$$
where $A$ is a sheaf on $\Spec(k)_\fppf$ and $B$ is a sheaf on $X_\fppf$. 
This is a particular case of the spectral sequence of composed functors,
namely $\Ga(X,-)$ and $\Hom_k(A,-)$, using that $\pi_X^*$
is a left adjoint to $\pi_{X*}$, and that $\pi_{X*}$ sends injective sheaves on $X_\fppf$
to injective sheaves on $\Spec(k)_\fppf$. 
The last property is a consequence of the fact that $\pi_X^*$
is exact, see \cite[Remark III.1.20]{EC} which refers to
\cite[Prop.~II.2.6]{EC}. See also \cite[\S 2.1.3]{CTS21} for a summary.

Since $\pi_X^*(\GG^\vee)=\GG^\vee_X$, we have the spectral sequence 
$$\Ext^p_{k}(\GG^\vee,R^q\pi_{X*}\G_{m,X})\Rightarrow 
\Ext^{p+q}_{X}(\GG^\vee_X,\G_{m,X}).$$
Since $X$ is proper, geometrically reduced and geometrically connected,
the natural morphism $\G_{m,k}\to \pi_{X*}\G_{m,X}$ is an isomorphism.
Thus the exact sequence of terms of low degree of our spectral sequence can be written as follows:
$$0\to \Ext^1_{k}(\GG^\vee,\G_{m,k})\to \Ext^1_{X}(\GG^\vee_X,\G_{m,X})\to
\Hom_{k}(\GG^\vee,{\bf Pic}_{X/k})$$
$$\to \Ext^2_{k}(\GG^\vee,\G_{m,k})\to \Ext^2_{X}(\GG^\vee_X,\G_{m,X}).$$
Using $x_0\in X(k)$ we obtain that the second and fifth arrows here are split injective.

We now consider the local-to-global spectral sequence of Ext-groups, see SGA 4, Exp.~V, (6.1.3):
$$\H^p(X,{\sE}xt^q_{X}(\GG^\vee_X,\G_{m,k}))\Rightarrow 
\Ext^{p+q}_{X}(\GG^\vee_X,\G_{m,X}).$$
By SGA 7, Exp.~VIII, Prop.~3.3.1, we have ${\sE}xt^1_{X}(\GG^\vee_X,\G_{m,k})=0$,
from which we obtain 
$$\Ext^1_{k}(\GG^\vee,\G_{m,k})\cong \H^1(k,\GG), \quad
\Ext^1_{X}(\GG^\vee_X,\G_{m,k})\cong \H^1(X,\GG).$$
Specialising at the base point $x_0$ we deduce the required isomorphism $\tau$. \hfill$\Box$

\medskip

It follows that if $p^n\GG=0$, then $\tau$ is an isomorphism $\H^1(X,\GG)_e\stackrel{\sim}\lra\Hom_k(\GG^\vee,S_X^\vee)$.

Let $S_X\otimes S_Y$ be the fppf sheaf of abelian groups on $\Spec(k)$ given by the tensor product of sheaves 
associated to the commutative group $k$-schemes $S_X$ and $S_Y$. 

\bthe \label{ku}
Let $X$ and $Y$ be pointed projective, geometrically reduced and 
geometrically connected varieties over a field $k$. Then there is an isomorphism
\begin{equation}
\Hom_k(S_X\otimes S_Y,\mu_{p^n})\cong\Hom_k(S_X, S_Y^\vee)\stackrel{\sim}\lra \H^2(X\times_k Y,\mu_{p^n})_\prim.
\label{eqq}
\end{equation}
\ethe
{\em Proof.} This follows from Proposition \ref{mur} and the natural isomorphisms
$$\H^1(X, S_Y^\vee)_e\cong\Hom_k(S_Y, S_X^\vee)\cong\Hom_k(S_X, S_Y^\vee)\cong
\Hom_k(S_X\otimes S_Y,\mu_{p^n}).$$
The first isomorphism is $\tau$ of Proposition \ref{tt} for $\GG= S_Y^\vee$. 
The second isomorphism is due to Cartier duality. The third isomorphism
is obtained by applying the functor of sections to the canonical isomorphism
$$\Hom(A,\Hom(B,C))\cong\Hom(A\otimes B,C)$$
in the category of fppf sheaves of abelian groups on $\Spec(k)$,
and noticing that $\Hom(S_Y,\mu_{p^n})\cong S_Y^\vee$ since $S_Y$ is annihilated by $p^n$.
\hfill$\Box$

\medskip

Following \cite[II (5.1.6)]{Ill} we define $\H^i(X,\Z_p(1))$ as $\varprojlim \H^i(X,\mu_{p^n})$ for $n\to\infty$.

\bco \label{si}
Let $X$ and $Y$ be pointed smooth, projective, geometrically integral varieties over a field $k$
of characteristic $p>0$. Then there is an isomorphism
\begin{equation}
\H^2(X\times_k Y,\Z_p(1))_\prim\cong\Hom_k(A[p^\infty],B^\vee[p^\infty]),
\label{fl}
\end{equation}
where $A[p^\infty]$ is the $p$-divisible group of the Albanese variety $A$ of $X$, and
$B^\vee[p^\infty]$ is the $p$-divisible group of the Picard variety $B^\vee$ of $Y$.
\eco
{\em Proof.} We have an exact sequence of group $k$-schemes
$$0\to {\rm{\bf Pic}}^0_{X/k}\to {\rm{\bf Pic}}_{X/k}\to {\rm{\bf NS}}_{X/k}\to 0,$$
which is the definition of ${\rm{\bf NS}}_{X/k}$, cf.~\cite[\S 5.1]{CTS21}.
The $k$-scheme ${\rm{\bf NS}}_{X/k}$ is \'etale, see SGA 3, IV$_A$, Prop.~5.5.1.
The group ${\rm{\bf NS}}_{X/k}(k^\s)={\rm{\bf NS}}_{X/k}(\bar k)=\NS(X_{\bar k})$ is finitely
generated by a theorem of N\'eron and Severi. Thus the cokernel of the map of 
group $k$-schemes ${\rm{\bf Pic}}^0_{X/k}[p^n]\to {\rm{\bf Pic}}_{X/k}[p^n]$
has bounded exponent.
Next, by Grothendieck [FGA6, \S 3], the Picard variety $A^\vee={\rm{\bf Pic}}^0_{X/k, {\rm red}}$
is a group subscheme of ${\rm{\bf Pic}}^0_{X/k}$ with finite cokernel, see 
\cite[\S 5.1.1]{CTS21}. We conclude that there is an exact sequence of finite commutative
group $k$-schemes
$$0\to A^\vee[p^n]\to S_X^\vee\to F_X\to 0,$$
where $F_X$ has bounded exponent. 
From the dual of this exact sequence and a similar sequence for $Y$ we obtain
the following exact sequence of abelian groups:
$$0\to \Hom_k(A[p^n],B^\vee[p^n])\to \Hom_k(S_X, S_Y^\vee)\to \Hom_k(S_X,F_Y)
\oplus\Hom_k(F_X^\vee,S_Y^\vee),$$
where homomorphisms are taken in the category of finite commutative $k$-groups.
We note that the last term in this sequence is annihilated by the maximum of the exponents
of $F_X$ and $F_Y$. This gives an isomorphism
$$\varprojlim \Hom_k(A[p^n],B^\vee[p^n])\cong \varprojlim \Hom_k(S_X, S_Y^\vee).$$
Thus passing to the projective limit in (\ref{eqq}) we obtain (\ref{fl}).
\hfill$\Box$

\medskip

We finish this section by interpreting the isomorphism (\ref{eqq})
of Theorem \ref{ku} in terms of certain canonical torsors on $X$ and $Y$.

For any $n\geq 1$ define a {\em universal $p^n$-torsor}\footnote{This notion
was introduced by Yang Cao in \cite{Cao20} (for the \'etale topology and for varieties without
a distinguished rational point),
inspired by universal torsors of Colliot-Th\'el\`ene and Sansuc \cite{CS}
and by calculations in \cite{SZ14}.} 
$\T_{X,p^n}\to X$ as an fppf
$X$-torsor with structure group $S_X$ and trivial fibre at $x_0$
such that the map $\tau$ from Proposition \ref{tt} sends 
the class $[\T_{X,p^n}]\in \H^1(X,S_X)_e$ to the natural injective map
$$S_X^\vee={\bf Pic}_{X/k}[p^n]\hookrightarrow {\bf Pic}_{X/k}.$$
It is clear that $\T_{X,p^n}$ is unique up to isomorphism.

The isomorphism (\ref{eqq})
in Theorem \ref{ku} can be made explicit in terms of $\T_{X,p^n}$ and $\T_{Y,p^n}$, as follows.
The cup-product pairing
$$\H^1(X\times Y,S_X)\times \H^1(X\times Y,S_Y)\to \H^2(X\times Y, S_X\otimes S_Y)$$
gives rise to the pairing
$$\H^1(X,S_X)\times \H^1(Y,S_Y)\to \H^2(X\times Y, S_X\otimes S_Y).$$
Let us denote by 
$$[\T_{X,p^n}]\boxtimes[\T_{X,p^n}]\in \H^2(X\times Y, S_X\otimes S_Y)_\prim$$
the value of the last pairing on the classes $[\T_{X,p^n}]$ and $[\T_{Y,p^n}]$.  
Define 
$$\e\colon \Hom_k(S_X\otimes S_Y,\mu_{p^n})\to \H^2(X\times_k Y,\mu_{p^n})_\prim$$
as the map sending a homomorphism $\psi\colon S_X\otimes S_Y\to\mu_{p^n}$
of sheaves on $\Spec(k)_\fppf$
to $\psi_*\big([\T_{X,p^n}]\boxtimes[\T_{X,p^n}]\big)$.

\bpr \label{last}
Let $X$ and $Y$ be pointed projective, geometrically reduced and 
geometrically connected varieties over 
a field $k$. The isomorphism $(\ref{eqq})$ is given by the map $\e$.
\epr
{\em Proof.} The second proof of \cite[Thm.~5.7.7 (ii)]{CTS21} on pp.~161--162 works
in our situation. We reproduce this argument for the convenience of the reader.

For a finite commutative group $k$-scheme $\GG$ such that $p^n\GG=0$ we have 
a commutative diagram of pairings:
$$\begin{array}{ccccc}
\H^1(X,\GG^\vee)_e&\times&\H^1(Y,\GG)_e&\to&\H^2(X\times_k Y,\mu_{p^n})\\
||&&\downarrow&&||\\
\H^1(X,\GG^\vee)_e&\times&\Ext^1_Y(\GG^\vee,\mu_{p^n})&\to&
\H^2(X\times_k Y,\mu_{p^n})\\
||&&||&&\uparrow\\
\H^1(X,\GG^\vee)_e&\times&\Ext^1_k(\GG^\vee,\tau_{\leq 1}{\bf R} {\pi_Y}_*\mu_{p^n})
&\to&\H^2(X,\tau_{\leq 1}{\bf R} {p_Y}_*\mu_{p^n})\\
||&&\downarrow&&\downarrow\\
\H^1(X,\GG^\vee)_e&\times&\Hom_k(\GG^\vee,S_Y^\vee)&\to&\H^1(X,S_Y^\vee)
\end{array}$$
The vertical map $\H^1(Y,\GG)\to \Ext^1_Y(\GG^\vee,\mu_{p^n})$
comes from the local-to-global spectral sequence (SGA 4, Exp.~V, (6.1.3))
$$\H^p(Y,{\sE}xt^q_Y(\GG^\vee,\mu_{p^n}))\Rightarrow\Ext^{p+q}_Y(\GG^\vee,\mu_{p^n}).$$
The first two pairings are compatible by \cite[Prop.~V.1.20]{EC}. 
The two lower pairings are natural, and the compatibility of the rest of the diagram is clear.
The composition of maps in the second column is the isomorphism $\tau$.

Since $Y$ is a pointed proper, geometrically reduced and geometrically connected variety over $k$, 
the object $\tau_{\leq 1}{\bf R} {p_Y}_*\mu_{p^n}$
of the bounded derived category of sheaves on $X_\fppf$ is the direct sum of $\mu_{p^n}$
in degree 0 and $S_Y^\vee$ in degree 1. Thus $\H^1(X,S_Y^\vee)$ is a direct summand of
$\H^2(X,\tau_{\leq 1}{\bf R} {p_Y}_*\mu_{p^n})$. Taking $\GG=S_Y$, 
the previous diagram gives rise to a commutative diagram of pairings
$$\begin{array}{ccccc}
\H^1(X,S_Y^\vee)_e&\times&\H^1(Y,S_Y)_e&\to&\H^2(X\times_k Y,\mu_{p^n})_\prim\\
||&&\tau\downarrow&&\uparrow\\
\H^1(X,S_Y^\vee)_e&\times&\Hom_k(S_Y^\vee,S_Y^\vee)&\to&\H^1(X,S_Y^\vee)_e
\end{array}$$
where both vertical arrows are isomorphisms of Propositions \ref{mur} and \ref{tt}.

Let $\psi\in\Hom_k(S_X\otimes S_Y,\mu_{p^n})$. 
Let $\varphi$ be the corresponding element in $\Hom(S_X,S_Y^\vee)$,
and let $\varphi^\vee\in\Hom(S_Y,S_X^\vee)$ be its dual. By construction,
the isomorphism $(\ref{eqq})$ sends $\psi$ to the image of $\tau^{-1}(\varphi^\vee)\in
\H^1(X,S_Y^\vee)_e$ in $\H^2(X\times_k Y,\mu_{p^n})_\prim$.
On the other hand, $\e(\psi)$ is the value of the top pairing of the last diagram on
$\varphi_*[\T_{X,p^n}]\in\H^1_\et(X,S_Y^\vee)_e$ and $[\T_{Y,p^n}]\in\H^1_\et(Y,S_Y)_e$.
Since $\tau([\T_{Y,p^n}])=\id\in\Hom(S_Y^\vee,S_Y^\vee)$, the commutativity of
the diagram shows that $\e(\psi)\in\H^2_\et(X\times_k Y,\Z/n)_\prim$
comes from $\varphi_*[\T_{X,p^n}]\in \H_\et^1(X,S_Y^\vee)$. Since
$\tau(\varphi_*[\T_{X,p^n}])$ is
the precomposition of $\tau([\T_X])=\id\in\Hom_k(S_X^\vee,S_X^\vee)$
with $\varphi^\vee\colon S_Y\to S_X^\vee$, we have
$\tau(\varphi_*[\T_{X,p^n}])=\varphi^\vee$. Thus $(\ref{eqq})$ coincides with
$\e$. \hfill $\Box$

\section{Brauer group of the product} \label{3}

For an abelian group $A$ the $p$-adic Tate module $T_p(A)$ is defined as the projective limit
$\varprojlim A[p^n]$ when $n\to\infty$. It is easy to see that $T_p(\Q_p/\Z_p)\cong\Z_p$
and that $T_p(M)=0$ if the abelian group $M$ has finite exponent.

\bthe \label{kanchi}
Let $X$ and $Y$ be pointed smooth, projective, geometrically integral varieties 
over a finitely generated field $k$ of characteristic $p>0$. Then we have the following statements.

{\rm (i)} The first Chern class gives an isomorphism
$$\Hom_k(A,B^\vee)\otimes\Z_p\stackrel{\sim}\lra\H^2(X\times_k Y,\Z_p(1))_\prim.$$

{\rm (ii)} We have $T_p(\Br(X\times_kY)_\prim)=0$.

{\rm (iii)} The abelian group $\Br(X\times_kY)\{p\}_\prim$ has finite exponent.
\ethe 
{\em Proof.} By a theorem of Chow (see \cite[Thm. 3.19]{Con06}), the natural map
$$\Hom_{k^\s}(A_{k^\s},B^\vee_{k^\s})\to\Hom_{\bar k}(A_{\bar k},B^\vee_{\bar k})$$
is an isomorphism. Hence we have natural isomorphisms:
$$
\Hom_k(A,B^\vee)\stackrel{\sim}\lra\Hom_{k^\s}(A_{k^\s},B^\vee_{k^\s})^\Ga
\stackrel{\sim}\lra\Hom_{\bar k}(A_{\bar k},B^\vee_{\bar k})^\Ga. 
$$
For a pointed projective, geometrically integral variety $(X,x_0)$ 
the natural map $\Pic(X)\to\Pic(X_{k^\s})^\Ga$ is an isomorphism \cite[Remark 5.4.3 (1)]{CTS21}.
Thus we obtain from \cite[Prop.~5.7.3]{CTS21} an isomorphism of abelian groups
$$
\Pic(X\times_kY)_\prim\cong \Hom_k(A,B^\vee). $$
Thus the primitive part of the Kummer exact sequence can be written as 
$$0\to \Hom_k(A,B^\vee)/p^n\stackrel{c_1}\lra \H^2(X\times_k Y,\mu_{p^n})_\prim\to 
\Br(X\times_kY)[p^n]_\prim\to 0.$$
The arrow marked $c_1$ is given by the first Chern class. Since $\Hom_k(A,B^\vee)$
is a finitely generated free abelian group, 
passing to the limit in $n$ and using Corollary \ref{si} we obtain an exact sequence
\begin{equation}
0\to \Hom_k(A,B^\vee)\otimes\Z_p\stackrel{c_1}\lra \Hom_k(A[p^\infty],B^\vee[p^\infty])\to 
T_p(\Br(X\times_kY)_\prim)\to 0.\label{mystery}
\end{equation}
De Jong's theorem (the crystalline Tate conjecture) \cite[Thm.~2.6]{dJ98} says
that the natural action of morphisms of abelian varieties on torsion points induces an
isomorphism
$$ \Hom_k(A,B^\vee)\otimes\Z_p\stackrel{\sim}\lra\Hom_k(A[p^\infty],B^\vee[p^\infty]).$$
This implies that the source and the target of the map $c_1$ are finitely generated
$\Z_p$-modules of the same rank. Since $T_p(\Br(X\times_kY)_\prim)$ is torsion-free,
the map $c_1$ must be an isomorphism, so $T_p(\Br(X\times_kY)_\prim)=0$. This proves (i) and (ii).

Let us prove (iii). For a finite extension $k'/k$
a standard restriction-corestriction argument \cite[Prop.~3.8.4]{CTS21}
 shows that the kernel of the natural map
$$\Br(X\times_kY)_\prim\to\Br(X_{k'}\times_{k'}Y_{k'})_\prim$$
is annihilated by $[k':k]$. Thus it is enough to prove (iii) after replacing $k$
by a finite field extension. In particular, we can assume that we have an isomorphism
$$\Hom_k(A,B^\vee)\stackrel{\sim}\lra\Hom_{\bar k}(A_{\bar k},B^\vee_{\bar k}).$$

Consider the commutative diagram with exact rows
$$\xymatrix{0\ar[r]&\Hom_{\bar k}(A_{\bar k},B^\vee_{\bar k})/p^n\ar[r]&
\H^2(X_{\bar k}\times_{\bar k}Y_{\bar k},\mu_{p^n})_\prim\ar[r]&
\Br(X_{\bar k}\times_{\bar k}Y_{\bar k})[p^n]_\prim\ar[r]&0\\
0\ar[r]&\Hom_k(A,B^\vee)/p^n\ar[r]\ar[u]^\cong&
\H^2(X\times_{k}Y,\mu_{p^n})_\prim\ar[r]\ar[u]&
\Br(X\times_{k}Y)[p^n]_\prim\ar[r]\ar[u]&0}$$
Comparing isomorphisms (\ref{eqq}) for $k$ and $\bar k$, 
we see that the middle vertical map is injective. Now the snake lemma
gives the injectivity of the right-hand map, hence
$\Br(X\times_kY)\{p\}_\prim$ is a subgroup of 
$\Br(X_{\bar k}\times_{\bar k}Y_{\bar k})\{p\}_\prim$.
By Theorem \ref{appendix: brauer finite} of the appendix, 
the group $\Br(X_{\bar k}\times_{\bar k}Y_{\bar k})\{p\}$
is the direct sum of an abelian $p$-group of  finite
exponent and finitely many copies of $\Q_p/\Z_p$, hence the same is true for 
$\Br(X\times_kY)\{p\}_\prim$. Thus (ii) implies (iii).  \hfill$\Box$

\bthe \label{wed}
Let $X$ and $Y$ be smooth, projective, geometrically integral varieties 
over a finitely generated field $k$. Then the cokernel of the natural map
$$\Br(X)\oplus\Br(Y)\to\Br(X\times_kY)$$
has finite exponent.
\ethe
{\em Proof.} Since $X$ and $Y$ are smooth,
there is a finite separable field extension $k\subset k'$ such that 
$X(k')\neq\emptyset$ and $Y(k')\neq\emptyset$. 
We have a commutative diagram with natural vertical maps, and horizontal maps given by
restriction or corestriction:
$$
\xymatrix{\Br(X)\oplus\Br(Y)\ar[r]^{\res \ \ }\ar[d]&\Br(X_{k'})\oplus\Br(Y_{k'})\ar[r]^{\rm cores}\ar[d]&
\Br(X)\oplus\Br(Y)\ar[d]\\
\Br(X\times_kY)\ar[r]^{\res \ \ }&\Br(X_{k'}\times_{k'}Y_{k'})\ar[r]^{\rm cores}&\Br(X\times_kY)}
$$
It is well known that ${\rm cores}\circ\res$ is multiplication by $[k':k]$, see
\cite[\S 3.8]{CTS21}.
In view of the direct sum decomposition (\ref{dec}), the cokernel of the middle vertical map
is isomorphic to $\Br(X_{k'}\times_{k'}Y_{k'})_\prim$.
By Theorem \ref{kanchi} (for the $p$-primary part) and \cite[Thm.~B]{SZ14} (for the prime-to-$p$
part), there is a positive integer that annihilates $\Br(X_{k'}\times_{k'}Y_{k'})_\prim$. 
The theorem follows from these facts and the commutativity of the diagram. \hfill $\Box$

\bco \label{wed1}
Let $X$ and $Y$ be smooth, projective, geometrically integral varieties 
over a finitely generated field $k$ of characteristic exponent $p$.
Then the cokernel of each of the following natural maps is a direct sum of a finite
group and a $p$-group of finite exponent:

{\rm (i)} 
$\Br(X_{k^\s})^\Ga\oplus\Br(Y_{k^\s})^\Ga\to\Br(X_{k^\s}\times_{k^\s}Y_{k^\s})^\Ga$;

{\rm (ii)} 
$\Br(X_{k^\s})^k\oplus\Br(Y_{k^\s})^k\to\Br(X_{k^\s}\times_{k^\s}Y_{k^\s})^k$;

{\rm (iii)} 
$\Br(X_{\bar k})^k\oplus\Br(Y_{\bar k})^k\to\Br(X_{\bar k}\times_{\bar k}Y_{\bar k})^k$.
\eco
{\em Proof.} (i) For every positive integer $n$ coprime to $p$
the group $\Br(X_{k^\s})[n]$ is finite, see, e.g., \cite[Cor.~5.2.8]{CTS21}.
Thus it remains to bound the exponent of the cokernel of the map in (i).
We have a commutative diagram
$$\xymatrix{\Br(X_{k^\s})^\Ga\oplus\Br(Y_{k^\s})^\Ga\ar[r]&
\Br(X_{k^\s}\times_{k^\s}Y_{k^\s})^\Ga\\
\Br(X)\oplus\Br(Y)\ar[r]\ar[u]&\Br(X\times_kY)\ar[u]}$$
By \cite[Thm.~5.4.12]{CTS21}, the cokernel of 
right-hand vertical map has finite exponent.
By Theorem \ref{wed} the cokernel of the lower horizontal map has finite exponent. 
Now (i) follows from the commutativity of the diagram.

(ii) As in (i), it is enough to prove that the cokernel has finite exponent. This immediately follows
from Theorem \ref{wed}. The same proof gives (iii). \hfill$\Box$

\medskip

The above results can be applied to varieties dominated by products of curves.
Here we content ourselves with the following statement.

\bco \label{chennai}
Let $k$ be a finitely generated field of characteristic $p>0$, and let
$d$ be a positive integer not divisible by $p$.
Let $X\subset\P^3_k$ be the surface given by $F(x_0,x_1)=G(x_2,x_3)$,
where $F$ and $G$ are homogeneous forms of degree $d$ without multiple roots.
Then $\big(\Br(X)/\Br_0(X)\big)(p')$ is finite and $\big(\Br(X)/\Br_0(X)\big)\{p\}$
has finite exponent. 
\eco
{\em Proof.} Let $C_1$ and $C_2$ be the plane curves given by
$y^d=F(x_0,x_1)$ and $z^d=G(x_2,x_3)$, respectively. 
Let $S_1\subset C_1$ be given by $y=0$ and let $S_2\subset C_2$ be given by $z=0$.
The rational map from $C_1\times_k C_2$ to $X$ sending $(x_0:x_1:y)\times(x_2:x_3:z)$
to $(zx_0:zx_1:yx_2:yx_3)$ is the composition of the following rational maps 
\cite[\S 1, Remark 1.10]{SK}:
\begin{itemize}
\item the inverse of the blow-up $Z\to C_1\times_k C_2$ of $S_1\times S_2\subset C_1\times C_2$; 

\item the quotient morphism $Z\to Z/\mu_d$, where $\mu_d$ acts diagonally on $y$ and $z$; 

\item the contraction $Z/\mu_d\to X$
of the images of the strict transforms of $S_1\times_k C_2$ and $C_1\times_k S_2$.
\end{itemize}
\noindent We note that $Z$ and $Z/\mu_d$ are non-singular \cite[Lemma 1.4]{SK},
and that $Z\to C_1\times_k C_2$ and $Z/\mu_d\to X$ are birational morphisms.
By the birational invariance of the Brauer group \cite[Corollary 6.2.11]{CTS21},
we have $\Br(C_1\times_k C_2)\cong \Br(Z)$ and $\Br(X)\cong\Br(Z/\mu_d)$.

Since $\Br(C_{1,k^\s})=0$ and $\Br(C_{2,k^\s})=0$ \cite[Theorem 5.6.1 (iv)]{CTS21},
Corollary \ref{wed1}~(ii) implies that $\Br(Z_{k^\s})^k$ has finite exponent. 
The kernel of $\Br(Z_{k^\s}/\mu_d)\to \Br(Z_{k^\s})$ is killed by $d$
\cite[Proposition 3.8.4]{CTS21}, so $\Br(X_{k^\s})^k$ also has finite exponent,
and thus is a direct sum of a finite group and a $p$-group of finite exponent. Since
$\Br_1(X)/\Br_0(X)$ is finite by \cite[Cor.~16.3.4]{CTS21} or \cite[\S 2.1]{GS22},
the statement follows. \hfill$\Box$

\section{Abelian varieties} \label{A}

The following lemma may be well known to the experts; we give a proof because we could not find it
in the literature.

\ble \label{lala}
Let $A$ be an abelian variety over an algebraically closed field $k$. Let $p$ be a prime, 
possibly equal to ${\rm char}(k)$. For any integer $m$ the endomorphism
$[m]\colon A\to A$ acts on $\H^2_\fppf(A,\mu_{p^n})$ as $m^2$ for any $n\geq 1$.
\ele
{\em Proof.} In the case $p\neq{\rm char}(k)$ we can replace fppf cohomology by \'etale
cohomology. Since $[m]$ acts on $\H^1_\et(A,\mu_{p^n})\cong A^\vee(k)[p^n]$
as $m$, it acts on 
$\H^2_\et(A,\mu_{p^n})\cong\wedge^2\H^1_\et(A,\mu_{p^n})(-1)$ as $m^2$.

Now let $p={\rm char}(k)$. 
Considering the map $[p^n]\colon \O_A^\times\to \O_A^\times$ in the fppf and \'etale topologies
gives rise to a canonical isomorphism \cite[(II.5.1.4)]{Ill}
$$\H^i_\fppf(A,\mu_{p^n})\cong\H^{i-1}_\et(A,\O_A^\times/O_A^{\times p^n}).$$
There is a map of \'etale sheaves of abelian groups 
$d\log\colon \O_A^\times/O_A^{\times p^n}\to W_n\Omega^1_X$,
see \cite[Prop.~I.3.23.2]{Ill}. By \cite[Thm.~II.1.4, (II.1.3.3)]{Ill} for each $i\geq 0$ we have a
canonical isomorphism $\H^i_\cris(A/W_n)\cong \H^i_\et(A,W_n\Omega_A^\bullet)$.
We claim that the resulting map
$$d\log\colon \H^1_\et(A,\O_A^\times/O_A^{\times p^n})\to \H^2_\cris(A/W_n)$$
is injective. We sketch the proof referring to \cite[Thm.~1.7]{YY} for details.

The case $n=1$ is stated in \cite[Remarque II.5.17 (a)]{Ill}. It is
a consequence of the following two facts:

\smallskip

(1) the map $\H^0(A,Z^1_A)\to \H^0(A,\Omega^1_A)$ is surjective, 
where $Z^1_A:=\Ker[d\colon\Omega^1_A\to\Omega^2_A]$;

(2) the map $\H^1(A,Z^1_A)\to \H^2_{\rm dR}(A/k)$ induced by the natural morphism of complexes
$Z^1_A[-1]\to \Omega^\bullet_A$ is injective.

\smallskip

\noindent Property (1) is true for any commutative group scheme $A$. Indeed, for invariant
vector fields $X$ and $Y$ and an invariant differential $\omega$ we have
$$d\omega(X,Y)=X\big(\omega(Y)\big)-Y\big(\omega(X)\big)+\omega([X,Y])=0,$$
because $\omega(Y)$ and $\omega(X)$ are in $k$, and 
the Lie algebra of $A$ is abelian.

The map in (2) factors as $\H^1(A,Z^1_A)\to {\mathbb H}^2(A,\Omega^{\geq 1}_A)\to
{\mathbb H}^2(A,\Omega^\bullet_A)=\H^2_{\rm dR}(A/k)$. The second arrow is injective
because for abelian varieties
the Hodge-de Rham spectral sequence degenerates at the first page, by 
a theorem of Oda \cite[Prop.~5.1]{Oda69}. 
The injectivity of the first map can be checked using \v{C}ech cohomology \cite[Prop.~5.6]{YY}.

The case of $n\geq 2$ follows by induction in $n$ from the following commutative diagram with
exact rows:
$$\xymatrix{0\ar[r]&\H^2_\fppf(A,\mu_{p^m})\ar[r]\ar[d]^{d\log}&
\H^2_\fppf(A,\mu_{p^{m+n}})\ar[r]\ar[d]^{d\log}&\H^2_\fppf(A,\mu_{p^n})\ar[d]^{d\log}\\
0\ar[r]&\H^2_\cris(A/W_m)\ar[r]&\H^2_\cris(A/W_{m+n})\ar[r]&\H^2_\cris(A/W_n)}$$
The zero in the top row is due to the natural isomorphism
$\H^1_\fppf(A,\mu_{p^n})\cong A^\vee(k)[p^n]$
and the surjectivity of multiplication by $p^m$ on $A^\vee(k)$. The zero in the bottom row
follows from the isomorphisms $\H^i_\cris(A/W_n)\cong \H^i_\cris(A/W)/p^n$
which are consequences of the fact 
that the groups $\H^i_\cris(A/W)$ are torsion-free $W$-modules.

A canonical isomorphism
$\H^i_\cris(X/W)\cong\wedge^i\H^1_\cris(X/W)$ shows that $[m]$ acts on $\H^i_\cris(X/W)$
as $m^i$. Thus the proposition follows from the claim. \hfill $\Box$

\medskip

We can use Theorem \ref{kanchi} to give a shorter proof of a result of D'Addezio 
\cite[Thm.~5.2]{D'A}.

\bthe \label{da}
Let $A$ be an abelian variety over a finitely generated field $k$ of characteristic exponent $p$.
Then $\Br(A_{\bar k})^k$ is a direct sum of a finite group and a $p$-group of finite exponent.
\ethe
{\em Proof.} Let $m\colon A\times_kA\to A$ be the group law of $A$. Define
$\delta\colon\Br(A)\to\Br(A\times A)$ as $m^*-\pi_1^*-\pi_2^*$. It is immediate to check that
$\delta(\Br(A)_e)\subset\Br(A\times A)_\prim$. By \cite[Lemma 2.1, Prop.~2.2]{inv}
we have an exact sequence
\begin{equation}
0\to \Br(A)_e\cap\Br_A(A)\to \Br(A)_e\stackrel{\delta}\lra\Br(A\times A)_\prim, \label{in}
\end{equation}
where $\Br_A(A)$ is the invariant Brauer group of $A$.
The group $\Br(A\times A)_\prim$ has finite exponent by Theorem \ref{wed}.
The image of $\Br_A(A)$ in $\Br(A_{\bar k})$ is contained in $\Br_A(A_{\bar k})$, but
$\Br_A(A_{\bar k})$ is annihilated by 2. Indeed, on the one hand,
by Lemma \ref{lala} and the Kummer exact sequence,
$[-1]^*$ acts on $\Br(A_{\bar k})$ trivially. On the other hand, 
$[-1]^*$ acts on $\Br_A(A_{\bar k})$ as $-1$, see \cite[Prop.~2.2]{inv}.
We conclude from (\ref{in}) that $\Br(A_{\bar k})^k$ has finite exponent.
It remains to use the finiteness of $\Br(A_{\bar k})[n]$ where $n$ is coprime to $p$, see
\cite[Cor.~5.2.8]{CTS21}. \hfill$\Box$

\brem{\rm Since the Picard scheme of an abelian variety is smooth, the natural map
$\Br(A_{k^\s})\to\Br(A_{\bar k})$ is injective \cite[Thm.~5.2.5 (ii)]{CTS21}, \cite[Cor.~3.4]{D'A},
thus $\Br(A_{k^\s})^k\cong\Br(A_{\bar k})^k$. By \cite[Thm.~5.4.12]{CTS21} we conclude from Theorem \ref{da} that
$\Br(A_{k^\s})^\Ga$ is a direct sum of a finite group and a $p$-group of finite exponent.
}
\erem

Using similar ideas, we can give a simplified proof of the flat version of the Tate conjecture
for divisors proved by D'Addezio in \cite[Thm.~5.1]{D'A}.

\bthe \label{dada}
Let $A$ be an abelian variety over a finitely generated field $k$ of characteristic $p>0$.
The image of $\H^2(A,\Z_p(1))$ in $\H^2(A_{\bar k},\Z_p(1))^\Ga$ is contained in the image
of the first Chern class map ${\rm c}_1\colon\NS(A_{\bar k})^\Ga\otimes\Z_p\to
\H^2(A_{\bar k},\Z_p(1))^\Ga$. After replacing $k$ with a finite separable extension,
the two images become equal.
\ethe
{\em Proof.} We continue to write $\delta=m^*-\pi_1^*-\pi_2^*$. We have a commutative diagram
$$\xymatrix{\H^2(A,\Z_p(1))_e\ar[r]\ar[d]^\delta&\H^2(A_{\bar k},\Z_p(1))^\Ga\ar[d]^\delta&
\NS(A_{\bar k})^\Ga\otimes\Z_p\ar[l]_{ \ \ {\rm c}_1}\ar[d]^\delta_\cong\\
\H^2(A\times_kA,\Z_p(1))_\prim^{\rm sym}\ar[r]
\ar@{..>}@(l,l)[rru]
&\H^2(A_{\bar k}\times_{\bar k}A_{\bar k},\Z_p(1))_\prim^{{\rm sym,}\, \Ga}&
\NS(A_{\bar k}\times_{\bar k}A_{\bar k})_\prim^{{\rm sym},\,\Ga}\otimes\Z_p\ar[l]_{ \ \ {\rm c}_1}
}$$
where the superscript `sym' stands for the elements fixed by the permutation of factors
in $A\times_kA$ and $A_{\bar k}\times_{\bar k}A_{\bar k}$. To prove the first statement it is
enough to construct the dotted line such that the resulting diagram is still commutative.

Theorem \ref{kanchi} (i) gives an isomorphism 
\begin{equation}
\Hom_k(A,A^\vee)\otimes\Z_p\stackrel{\sim}\lra \Pic(A\times_kA)_\prim\otimes\Z_p
\stackrel{\sim}\lra \H^2(A\times_kA,\Z_p(1))_\prim.\label{odin}
\end{equation}
Here the first arrow sends $f\in \Hom_k(A,A^\vee)$ to $(\id,f)^*{\mathcal P}$,
where ${\mathcal P}$ is the Poincar\'e line bundle on $A\times_kA^\vee$. 
The second arrow is the first Chern class ${\rm c}_1$. If $f=f^\vee$, then
the image of $f$ lands in the symmetric subgroup of $\H^2(A\times_kA,\Z_p(1))_\prim$.
The same construction over $\bar k$ gives an isomorphism of $\Ga$-modules
$$ \Hom(A_{\bar k},A^\vee_{\bar k})\otimes\Z_p\stackrel{\sim}\lra 
\NS(A_{\bar k}\times_{\bar k}A_{\bar k})_\prim\otimes\Z_p\cong\Pic(A_{\bar k}\times_{\bar k}A_{\bar k})_\prim\otimes\Z_p,$$
which is clearly compatible with the first map of (\ref{odin}) and which gives this map
after restricting to the $\Ga$-invariant subgroups. We finally note that the isomorphism
of $\Ga$-modules $\Hom(A_{\bar k},A^\vee_{\bar k})^{\rm sym}\cong
\NS(A_{\bar k}\times_{\bar k}A_{\bar k})_\prim^{\rm sym}$
identifies the map of $\Ga$-modules
$$\delta\colon\NS(A_{\bar k})\to \NS(A_{\bar k}\times_{\bar k}A_{\bar k})_\prim^{\rm sym}$$
with the isomorphism 
$\NS(A_{\bar k})\stackrel{\sim}\lra\Hom(A_{\bar k},A^\vee_{\bar k})^{\rm sym}$
sending $L$ to $\varphi_L$. (This follows from
$(\id,\varphi_L)^*{\mathcal P}=m^*L\otimes\pi_1^*L^{-1}\otimes\pi_2^*L^{-1}$,
see \cite[Ch.~8]{Mum}, cf.~\cite[Prop.~6.1]{inv}.) Putting all of this together gives rise
to a dotted line in the diagram, which is the identity map on $\Hom_k(A,A^\vee)^{\rm sym}\otimes\Z_p$ once the source and the target are identified with this group. The resulting
diagram commutes.
This follows from the injectivity of the middle vertical map $\delta$: arguing as in the proof of the previous theorem, one shows that the kernel of this map is annihilated by 2, but $\H^2(A_{\bar k},\Z_p(1))$ is torsion-free \cite[Thm.~II.5.14] {Ill}.

Since $\NS(A_{\bar k})$ is finitely generated,
replacing $k$ by a finite separable extension we can ensure that the map 
$\Pic(A)\to\NS(A_{\bar k})^\Ga$ is surjective. Then the image
of the first Chern class map ${\rm c}_1\colon\NS(A_{\bar k})^\Ga\otimes\Z_p\to
\H^2(A_{\bar k},\Z_p(1))^\Ga$ is contained in the image of 
$\H^2(A,\Z_p(1))\to\H^2(A_{\bar k},\Z_p(1))^\Ga$. The second statement follows. \hfill$\Box$

\section{Kummer surfaces} \label{K}

Recall that the Picard scheme of a K3 surface is smooth, hence the natural map
$\Br(X_{k^\s})\to\Br(X_{\bar k})$ is injective \cite[Thm.~5.2.5 (ii)]{CTS21}, \cite[Cor.~3.4]{D'A}. 
This implies $\Br(X_{k^\s})^k\cong\Br(X_{\bar k})^k$.

\bpr \label{yi}
Let $k$ be a field of characteristic exponent $p\neq 2$.
Let $A$ be an abelian surface and let $X=\Kum(A)$ be the associated Kummer surface.
Then there is a natural isomorphism of $\Ga$-modules 
$\Br(X_{\bar k})\stackrel{\sim}\lra\Br(A_{\bar k})$.
\epr
{\em Proof.} This is proved in \cite[Prop.~2.1.1]{LZ}.
For all primes $\ell\neq p$ (including $\ell=2$) the proof 
of \cite[Prop.~1.3]{SZ12} shows that
$\Br(X_{\bar k})\{\ell\}\to\Br(A_{\bar k})\{\ell\}$ is an isomorphism.
In fact, for any $\ell\neq 2$ (including $\ell=p$ if $p>1$) the map 
$\Br(X_{\bar k})\{\ell\}\to\Br(A_{\bar k})\{\ell\}$
is injective with image $\Br(A_{\bar k})\{\ell\}^{[-1]^*}$ by 
\cite[Thm.~3.8.5]{CTS21}. In view of the Kummer sequence, it suffices to show
that $[-1]$ acts on $\H^2_\fppf(A_{\bar k},\mu_{\ell^n})$ trivially.
This was proved in Lemma \ref{lala}. \hfill $\Box$

\bco \label{er}
Let $k$ be a field of characteristic exponent $p\neq 2$.
Let $A$ be an abelian surface and let $X=\Kum(Y)$ be the associated Kummer surface.
For all odd primes $\ell$ (including $\ell=p$ if $p>1$) there are natural isomorphisms of abelian groups 
$$\Br(X_{\bar k})\{\ell\}^k\stackrel{\sim}\lra\Br(A_{\bar k})\{\ell\}^k.$$
\eco
{\em Proof.} This is proved in \cite[Thm.~2.4]{SZ12}. 
Let us give this proof for the convenience of the reader.
By Proposition \ref{yi}, the map is injective, so it remains to prove that it is surjective.
For any odd $\ell$ we have a direct sum decomposition
$$\Br(A)\{\ell\}=\Br(A)\{\ell\}^+\oplus \Br(A)\{\ell\}^-,$$
where  $\Br(A)\{\ell\}^+=\Br(A)\{\ell\}^{[-1]^*}$ is the $[-1]^*$-invariant
subgroup and $\Br(A)\{\ell\}^-$ is the $[-1]^*$-antiinvariant
subgroup. By the proof of Proposition \ref{yi}, the action of $[-1]$ on $\Br(A_{\bar k})$ is trivial,
thus any element of $\Br(A_{\bar k})\{\ell\}^k$ lifts to an element of $\Br(A)\{\ell\}^+$.
The last group is the image of $\Br(X)\{\ell\}$ by \cite[Thm.~3.8.5]{CTS21}. \hfill $\Box$

\bco \label{san}
Let $k$ be a finitely generated field of characteristic exponent $p\neq 2$.
Let $A$ be an abelian surface and let $X=\Kum(A)$ be the associated Kummer surface.
Then each of the groups $\Br(X_{k^\s})^\Ga$ and $\Br(X_{k^\s})^k\cong\Br(X_{\bar k})^k$ 
is a direct sum of a finite group and and a $p$-group of finite exponent.
\eco
{\em Proof.} For any K3 surface $X$ over $k$, 
the finiteness of $\Br(X_{k^\s})(p')^\Ga$ and $\Br(X_{k^\s})(p')^k$
was proved in \cite[Thm.~1.2]{SZ08} for $p=1$ and in \cite[Thm.~1.3]{SZ15} for $p>2$.
For $p>2$ the statements for $p$-primary torsion follow from 
Theorem \ref{da} and Corollary \ref{er}, using \cite[Thm.~5.4.12]{CTS21}
which says that the quotient of $\Br(X_{k^\s})^\Ga$ by $\Br(X_{k^\s})^k$ 
is a direct sum of a finite group and and a $p$-group of finite exponent.
\hfill $\Box$

\medskip

\noindent{\bf Example 1}.
The group $\Br(A_{k^\s})[p]^k$ may well be infinite \cite[Cor.~5.4]{D'A}.
Let $E$ be a supersingular elliptic curve over an {\em infinite} finitely 
generated field $k$ of characteristic $p>0$, and let $A=E\times_kE$.
The group scheme $E[p]$ is an extension of $\alpha_p$ by $\alpha_p$, hence
there is an injective map of abelian groups $\End_k(\alpha_p)\to \End_k(E[p])$ which sends an
endomorphism $\phi\colon\alpha_p\to\alpha_p$ to the composition
$$E[p]\to\alpha_p\stackrel{\phi}\lra\alpha_p\to E[p].$$
By Theorem \ref{ku} we have $\H^2(A,\mu_p)_\prim\cong\End_k(E[p])$, hence
$$\Br(A)[p]_\prim\cong \End_k(E[p])/\big(\End_k(E)/p\big).$$
Since $\End_k(\alpha_p)\cong k$, we have compatible
 injective homomorphisms $k\to \End_k(E[p])$ and $\bar k\to \End_{\bar k}(E[p])$.
Since $\End_{\bar k}(E)/p$ is finite, the image of $\Br(A)[p]$ in 
$$\Br(A_{\bar k})[p]_\prim\cong \End_{\bar k}(E[p])/\big(\End_{\bar k}(E)/p\big)$$
is infinite.
Now let $p\neq 2$. Then we can consider the Kummer surface $X=\Kum(A)$ over $k$.
Corollary \ref{er} implies that $\Br(X_{k^\s})[p]^k\cong\Br(X_{\bar k})[p]^k$
is infinite. This gives an example of a 
K3 surface with an infinite transcendental Brauer group, answering \cite[Questions 1, 2]{SZ08}
in the negative.

\medskip

\noindent{\bf Example 2}.
D'Addezio also shows that in the case of finite characteristic, the group 
$\Br(A_{\bar k})^\Ga$ does not always have finite exponent \cite[Cor.~6.7]{D'A}. Take
$A=E\times_k E$, where $E$ is an  elliptic curve over $k$ 
whose $j$-invariant is transcendental over $\F_p$, so that $E$ is ordinary and 
$\End_{\bar k}(E)\cong\Z$. By (\ref{mystery}),
the quotient of $\End(E_{\bar k}[p^\infty])$
by $\End(E_{\bar k})\otimes\Z_p\cong\Z_p$ is contained in $T_p(\Br(A_{\bar k}))$. 
Thus $T_p(\Br(A_{\bar k}))^\Ga$ contains the quotient of $\End(E_{\bar k}[p^\infty])^\Ga$
by $\Z_p$, so it is enough to show that the rank of the $\Z_p$-module
$\End(E_{\bar k}[p^\infty])^\Ga$ is at least 2. Since $E$ is ordinary,
the $p$-divisible group $E_{\bar k}[p^\infty]$ has two slopes: 0 and 1.
Let $k^{\rm perf}$ be the perfect closure of $k$ in $\bar k$. By the splitness of the 
connected-\'etale sequence over perfect fields, the $p$-divisible group
$E_{k^{\rm perf}}[p^\infty]$ 
is isomorphic to the direct sum of its connected and \'etale parts.
It follows that $\End(E_{\bar k}[p^\infty])^\Ga\cong \End(E_{k^{\rm perf}}[p^\infty])$ 
contains $\Z_p^{\oplus 2}$.

As before, if $p\neq 2$, then for $X=\Kum(A)$ we obtain from Proposition \ref{yi} that
$\Br(X_{\bar k})^\Ga$ does not have finite exponent.

\appendix

\renewcommand*{\thefootnote}{\fnsymbol{footnote}}
\section{Appendix, by Alexander Petrov\protect\footnote{A.P. was supported by the Clay Research Fellowship.}}
\renewcommand*{\thefootnote}{\arabic{footnote}}

Let $k$ be an algebraically closed field of characteristic $p>0$.
We write $W=W(k)$ for the ring of Witt vectors of $k$ and $K$ for the field of fractions of $W$.

For a smooth proper variety $X$ over $k$ we denote by 
$\rho=\dim_\Q(\NS(X)\otimes\Q)$ the Picard number of $X$. For $i\geq 0$, let $r_i$ be the dimension of the $\Q_p$-vector space 
$(\H^i_\cris(X/W)\otimes K)^{F=p}$.

Consider the complex of weight $1$ syntomic cohomology of $X$: 
$$
R\Ga(X,\Z_p(1)):=R\varprojlim\big(R\Gamma_{\fppf}(X,\mu_{p^n})\big)
$$
Here $R\varprojlim$ is the derived inverse limit \cite[08TC]{stacks-project} of the system of objects $R\Gamma_{\fppf}(X,\mu_{p^n})\in D(\Z_p)$ of the derived category of $\Z_p$-modules. In fact, each individual cohomology group $\H^i(X,\Z_p(1)):=\H^i(R\Gamma(X,\Z_p(1)))$ is isomorphic to  $\varprojlim \H^i_{\fppf}(X,\mu_{p^n})$ by the proof of \cite[Thm.~II.5.5]{Ill}.

The following result is well known to the experts, and can be deduced from
\cite[Prop.~II.5.9]{Ill} and its proof. For the convenience of the reader we give a
complete proof below.

\bthe\label{appendix: brauer finite}
Let $X$ be a smooth proper variety over $k$. 
Then there is an isomorphism of abelian groups
$$\Br(X)\{p\}\cong(\Q_p/\Z_p)^{\oplus (r_2-\rho)}\oplus \H^3(X,\Z_p(1))\{p\},$$
where the group $\H^3(X,\Z_p(1))\{p\}$ is annihilated by a power of $p$.
\ethe

Syntomic cohomology modules $\H^i(X,\Z_p(1))$ are finitely generated over $\Z_p$ 
for $i\leq 2$, see \cite[Prop.~II.5.9]{Ill}, but not always for $i\geq 3$ 
(cf.~the example of a supersingular K3 surface in 
\cite[II.7.2]{Ill}). However, they satisfy a weaker finiteness property that we will use to deduce Theorem \ref{appendix: brauer finite}:

\ble[Illusie--Raynaud]\label{appendix: syntomic torsion finite}
For each $i\geq 0$, the $\Z_p$-module $\H^i(X,\Z_p(1))$ 
is isomorphic to $\Z_p^{\oplus r_i}\oplus \H^i(X,\Z_p(1))\{p\}$, and the module
$\H^i(X,\Z_p(1))\{p\}$ is annihilated by a power of $p$.
\ele
{\em Proof.} This result follows from 
Th\'eor\`eme IV.3.3 (b) and Corollaire IV.3.6 of \cite{Ill-Ray}. Here we give a self-contained argument that uses only the more basic properties of the de Rham--Witt complex. The statement is clear for $i=0$, so we assume $i\geq 1$.

By \cite[Thm.~II.5.5]{Ill}, syntomic cohomology groups fit into the long exact sequence
\begin{multline}\label{appendix: syntomic to dRW sequence}
\ldots\to \H^{i-1}_{\Zar}(X,W\Omega^1_X)\xrightarrow{1-F}\H^{i-1}_{\Zar}(X,W\Omega^1_X)\to \H^{i+1}(X,\Z_p(1)) 
\to \\ \to\H^{i}_{\Zar}(X,W\Omega^1_X)\xrightarrow{1-F}\H^{i}_{\Zar}(X,W\Omega^1_X)\to\ldots
\end{multline}
where $W\Omega^1_X$ is the sheaf of de Rham--Witt differential forms, and $F:W\Omega^1_X\to W\Omega^1_X$ is its semi-linear Frobenius endomorphism \cite[I.2.E]{Ill}.

While $\H^j_{\Zar}(X,W\Omega^i_X)$ is not always finitely generated as a $W$-module, it is isomorphic to a direct sum of a finitely generated free $W$-module and a module annihilated by a power of $p$, by \cite[Thm.~II.2.13]{Ill}.

If $M$ is a finite free $W$-module equipped with a Frobenius-linear endomorphism $F:M\to M$, then $1-F:M\to M$ is surjective by \cite[Lemme II.5.3]{Ill} (this is the only place where we use that $k$ is algebraically closed rather than just perfect). The kernel $M^{F=1}:=\ker(1-F:M\to M)$ is a finite free $\Z_p$-module (cf. \cite[Lemme II.5.11]{Ill}) because the natural map $M^{F=1}\otimes_{\Z_p}W\to M$ is an injection, which follows from the fact that $W^{F=1}$ equals $\Z_p$ . 

More generally, this implies that if $\widetilde{M}$ is a $W$-module isomorphic to a direct sum of a finite free module and a module annihilated by a power of $p$, then for a Frobenius-linear endomorphism $F:\widetilde{M}\to \widetilde{M}$ the cokernel of $1-F$ is annihilated by a power of $p$, and its kernel is a direct sum of a finite free $\Z_p$-module and a $\Z_p$-module annihilated by a power of $p$.

Therefore the sequence (\ref{appendix: syntomic to dRW sequence}) gives that $\H^i(X,\Z_p(1))$ fits into a short exact sequence
\begin{equation}
0\to T'\to \H^i(X,\Z_p(1))\to T\oplus\Z_p^{\oplus r}\to 0
\end{equation}for some integer $r$, where both $T$ and $T'$ are annihilated by powers of $p$.
This implies that $\H^i(X,\Z_p(1))$ is isomorphic to $\Z_p^{\oplus r}\oplus \H^i(X,\Z_p(1))\{p\}$ with $\H^i(X,\Z_p(1))\{p\}$ annihilated by a power of $p$.

By \cite[Thm.~II.5.5 (5.5.3)]{Ill} we have an isomorphism of $\Q_p$-vector spaces
$$\H^i(X,\Z_p(1))\otimes\Q_p\cong(\H^i_\cris(X/W)\otimes K)^{F=p},$$
thus $r=r_i$. \hfill $\Box$

\brem {\rm
Illusie and Raynaud prove in \cite[Thm.~IV.3.3]{Ill-Ray} that the maps $1-F$ in the long exact sequence (\ref{appendix: syntomic to dRW sequence}) are in fact surjective, but we did not use this fact in the above proof.
}\erem

\noindent{\em Proof of Theorem} \ref{appendix: brauer finite}
For each $n$, we have a distinguished triangle in the derived category of $\Z_p$-modules
\begin{equation}\label{appendix: univ coefficients triangle}
R\Gamma(X, \Z_p(1))\xrightarrow{p^n}R\Gamma(X, \Z_p(1))\to R\Gamma_{\fppf}(X,\mu_{p^n})
\end{equation}
obtained from the distinguished triangles
$$
R\Gamma_{\fppf}(X,\mu_{p^m})\xrightarrow{p^n} R\Gamma_{\fppf}(X,\mu_{p^{n+m}})\xrightarrow{} R\Gamma_{\fppf}(X,\mu_{p^n})
$$
by passing to the inverse limit over all $m$. For all $i,n$ the triangle (\ref{appendix: univ coefficients triangle}) induces the short exact sequences
\begin{equation}
0\to \H^i(X,\Z_p(1))/p^n\to \H^i(X, \mu_{p^n})\to \H^{i+1}(X,\Z_p(1))[p^n]\to 0.
\end{equation}
For each $i$, passing to the direct limit along the maps induced by $\mu_{p^n}\hookrightarrow\mu_{p^{n+1}}$ we get the short exact sequence
\begin{equation}
0\to \H^i(X,\Z_p(1))\otimes_{\Z_p}\Q_p/\Z_p\to \varinjlim\H^i(X,\mu_{p^n})\to \H^{i+1}(X,\Z_p(1))\{p\}\to 0.
\end{equation}
By Lemma \ref{appendix: syntomic torsion finite} the group $\H^i(X,\Z_p(1))\otimes_{\Z_p}\Q_p/\Z_p$ is isomorphic to $(\Q_p/\Z_p)^{\oplus r_i}$, and the group $\H^{i+1}(X,\Z_p(1))\{p\}$ is annihilated by a power of $p$. The abelian group $\Q_p/\Z_p$ is divisible, hence injective,
thus $\varinjlim\H^i(X,\mu_{p^n})$ is isomorphic to the direct sum of 
$(\Q_p/\Z_p)^{\oplus r_i}$ and $\H^{i+1}(X,\Z_p(1))\{p\}$.

On the other hand, for all $n$ we have short exact sequences
\begin{equation}
0\to \Pic(X)/p^n\to \H^2(X,\mu_{p^n})\to \Br(X)[p^n]\to 0
\end{equation}
which induce, after passing to the direct limit, a surjection $
\varinjlim\H^2(X,\mu_{p^n})\to \Br(X)\{p\}$ with kernel
 $\varinjlim\Pic(X)/p^n\cong(\Q_p/\Z_p)^\rho$. 
Hence $$\Br(X)\{p\}\cong\mathrm{Coker}(\alpha)\oplus \H^3(X, \Z_p(1))\{p\}$$ for some injection $\alpha\colon(\Q_p/\Z_p)^{\oplus \rho}\hookrightarrow (\Q_p/\Z_p)^{\oplus r_2}$, which proves the theorem. \hfill $\Box$
\bigskip

MIT, Department of Mathematics, 77 Massachusetts Avenue, Cambridge, MA 02139 USA

\texttt{alexander.petrov.57@gmail.com}

 Department of Mathematics, 
South Kensington Campus,
Imperial College London,
SW7~2AZ United Kingdom \  \ and \  \ 
Institute for the Information Transmission Problems,
Russian Academy of Sciences,
Moscow, 127994 Russia

\texttt{a.skorobogatov@imperial.ac.uk}

\end{document}